\documentclass{amsart}
\usepackage{amssymb}
\usepackage{amsfonts}
\usepackage{hyperref}

\theoremstyle{plain}

\numberwithin{equation}{section}

\input{tcilatex}

\newcommand{\R}{{\mathbb{R}}}
\newcommand{\N}{{\mathbb{N}}}
\begin{document}
\title[Pendants to the Euler Beta function]{Pendants to the Euler Beta function\\
}
\author{Martin Himmel}
\curraddr{%Freie Werkschule Meißen,
 %Zscheilaer Str. 19, 
 %01662 Meißen, 
 Germany}
\email{martin.himmel@gmail.com}

\begin{abstract}
Motivated by the integral representation of the Euler Beta function,
we introduce its Cauchy siblings and investigate some of their properties.
Two of these newly introduced functions happen to coincide with some classical means,
such as the arithmetic or the logarithmic Cauchy one.
Although the bivariable generalizations of Beta functions are obtained by elemantary integration,
it seems difficult to obtain closed formulas for more than two variables.
The questions whether these Cauchy Beta functions belong to their respective class of Cauchy quotients
is addressed and answered positively in the case of the Euler Beta function,
but postponed to a future paper for all the other cases.  
\end{abstract}

\maketitle

\QTP{Body Math}
$\bigskip $

\section{Introduction}
\begin{quote}
	First children are born,\\
	later they perform miracles.
	
	unknown Author
\end{quote}
A couple of years ago the relationsip between the Beta function 
and the Euler gamma function 
\begin{equation*}
\mathcal{B}(x,y)=\frac{\Gamma(x) \Gamma(y)}{\Gamma(x+y)}	
\end{equation*}
gave rise to introduce a class called beta-type functions,
which are functions of the form
\begin{equation*}
	{B}_f(x,y)=\frac{f(x) f(y)}{f(x+y)},	
\end{equation*}
where $f$ is a positive function defined on a suitable real interval.
In this connection it was natural to ask for means in the class of beta-type functions, 
which were determined in \cite{MatHim2} and \cite{MatHim4}.
The form of beta-type functions suggested very strongly to introduce three more classes of functions,
each of them related to one of the four Cauchy functional equations (cf. \cite{Kuczma2}).\\
It turned out that in three of these classes, namely in beta-type functions (also called exponential Cauchy quotients), in logarithmic and in multiplicative Cauchy quotients there exists a unique mean of $k$ variables where $k$ is an arbitrary positive integer (cf. \cite{MatHim5}),
but in the class of additive Cauchy quotients there is no mean for any number of variables greater or equal than two.
The classes of means which appeared naturally in this area of research even led to unexpected applications in stochastic analysis (cf. \cite{Burai}).   

In this paper we take a closer look to the integral representation of the Beta function.
In a nutshell, we observe 
that the integrand of the Beta function consists basically of two smooth exponential functions
where the bases satisfy a certain duality relation.
To avoid any misunderstanding, we go into medias res.

The Euler Beta function $\mathcal{B}:\left(
0,+\infty \right) ^{2}\rightarrow \left( 0,+\infty \right)$ is defined by 
\begin{equation}
\mathcal{B}\left( x,y\right) =
\int_0^1 {t^{x-1} (1-t)^{y-1}\, dt},\text{ \ \ \ \ \ }x,y \in (0, +\infty).
\label{eq:Beta}
\end{equation}%
Intuitively, the ingredients of the Beta function are two smooth families of exponential functions of base $t$ and $1-t$, respectively, which are 'averaged together' in the sense that they are multiplied and then integrated over the unit interval.
More precisely, for fixed $t \in [0,1]$,
the exponential function $\R \ni x \mapsto t^{x-1} \in (0,+\infty)$ is multiplied by its dual exponential function
$\R \ni y \mapsto (1-t)^{y-1} \in (0,+\infty)$ (\textit{dual} in the sense that their respective bases add to one)
%both restricted to $(0,+\infty)$,
and integrated over the interval $[0,1]$.

In a similar manner we introduce in this paper about three more Beta functions, each of them related to one of the other Cauchy functional equations.

\section{Multiplicative Beta function}
Since the ingredients of the Beta function, i.e. the functions under the integral in \eqref{eq:Beta}, are smooth exponential functions,
it seems natural to consider the functions
$(1,  +\infty) \ni x \mapsto (x-1)^t \in (0,  +\infty)$ and $(1, +\infty) \ni y \mapsto (y-1)^{1-t} \in (0,  +\infty)$,
two families of smooth multiplicative functions dual to each other in the sense that the sum of their exponents is one, 
as the  ingredients of the multiplicative Beta function 
%$\mathcal{B}_m:\left(1,+\infty \right) ^{2}\rightarrow \left( 0,+\infty \right)$ 
$\mathcal{M}:\left(1,+\infty \right) ^{2}\rightarrow \R$
defined by

\begin{equation}
%\mathcal{B}_m\left( x,y\right) :=
\mathcal{M}\left( x,y\right)=
\int_0^1 {(x-1)^t (y-1)^{1-t}\, dt}, \text{ \ \ \ \ \ }x,y \in (1, +\infty).
\label{eq:Beta_m:def}
\end{equation}
Unlike to the case of the Euler Beta function, which we should call exponential Beta function from a more systematic point of view,
the multiplicative Beta function $\mathcal{M}$ can be expressed by elementary functions, since, for all $x,y \in (1, +\infty), x \neq y$, we have
\begin{eqnarray}
%\mathcal{B}_m
\label{eq:B_m}
\mathcal{M}\left( x,y\right) & = & (y-1) \int_0^1 {\left(\frac{x-1}{y-1}\right)^t\, dt}\\ \nonumber
& = &  (y-1) \left[  \frac{\left(\frac{x-1}{y-1}\right)^t}{\log{\left(\frac{x-1}{y-1}\right)}} \right]_0^1 \\ \nonumber
& = &\frac{x-y}{ \log{\left(\frac{x-1}{y-1}\right)}},
\end{eqnarray}
and $\mathcal{M}\left( x,x\right) = x-1$ for $x>1$ implying that $\mathcal{M}$ is the logarithmic mean where the arguments are both reduced by $1$.

\section{Additive Beta function}
The building blocks of the additive Beta function are
the two families 
$\R \ni x \mapsto t (x-1) \in \R$ and $\R  \ni y \mapsto ({1-t})(y-1)  \in \R$.
As before, we could argue that these $t$-families of continuous additive functions should be multiplied and then integrated to obtain
the additive Beta function
$\mathcal{A}_{1}: \R^2 \to \R$ defined by
\begin{equation*}
%\mathcal{B}_{a1} \left( x,y\right):=
\mathcal{A}_{1} \left( x,y\right):=\int_0^1 {t (x-1) \cdot ({1-t})(y-1) \, dt},
%\label{eq:}
\end{equation*}
which, for all $x,y \in\R$, simplifies to

\begin{eqnarray*}
%\mathcal{B}_{a1} \left( x,y\right)
\mathcal{A}_{1} \left( x,y\right)
&=&(x-1)(y-1) \int_0^1 {t(1-t)\, dt}\\
&=&(x-1)(y-1) {\left[ \frac{1}{2} t^2-\frac{1}{3} t^3\right]_0^1}\\
&=&\frac{1}{6} (x-1)(y-1).
%\label{eq:}
\end{eqnarray*}
But one could object that in the construction of the preceding Beta functions the motivation to multiply the building blocks, 
i.e. the two families of exponential functions (multiplicative functions) in the case of
the Euler Beta function (in the case of multiplicative Beta function, respectively) stems from the fact that 
in the exponential Cauchy equation 
\begin{equation*}
f(x+y)=f(x) \cdot f(y), \qquad x,y \in \R,
%\label{eq:}
\end{equation*}
(in the case of the multiplicative Cauchy equation $f(xy)=f(x)f(y)$ for all $x,y \in (0, +\infty)$, respectively)
multiplication of the function $f$ occurs on the right hand side.

If this is the right explanation
why the respective ingredients of exponential and multiplicative Beta functions are multiplied under the integral sign,
we should rather add the building blocks in the case of the additive Beta function,
since for an additive function we have
\begin{equation*}
g(x+y)=g(x) + g(y), \quad x,y \in \R,
%\label{eq:}
\end{equation*}
and addition occurs on the right hand side.
Thus, we may rather consider $\mathcal{A}_{2}: \R^2 \to \R$ defined by
\begin{equation}
%\mathcal{B}_{a2} \left( x,y\right):=
\mathcal{A}_{2} \left( x,y\right)=\int_0^1 {\left[t (x-1) + ({1-t})(y-1)\right] \, dt},
%\label{eq:}
\end{equation}
as the natural additive counterpart of the Euler Beta function.
Surprisingly, the function $\mathcal{A}_{2}$ coincides with the arithmetic mean,
since, for all $x,y \in \R$, it holds
\begin{eqnarray*}
%\mathcal{B}_{a2} \left( x,y\right)=
\mathcal{A}_{2}\left( x,y\right)&=& \left[{\frac{t^2}{2} (x-1) + \left({t-\frac{t^2}{2}}\right)(y-1)}\right]_0^1\\
&=&{\frac{1^2}{2} (x-1) + \left({1-\frac{1^2}{2}}\right)(y-1)}\\
&=&\frac{1}{2} (x-1) + \frac{1}{2} (y-1)\\
&=&A(x-1,y-1),
%\label{eq:}
\end{eqnarray*}
where $A:\R^2 \to \R$, defined by $A(x,y)=\frac{x+y}{2}$ for all $x,y \in \R$, is the arithmetic mean.

\section{Logarithmic Beta function}
Similarly as before, we consider the $t$-families $x \mapsto t \log{x}$ and $y \mapsto (1-t)\log{y}$ as the building blocks 
of the logarithmic Beta function $\mathcal{L}_{1}: (1,+\infty) \to \R$ defined by

\begin{equation}
\mathcal{L}_{1} \left( x,y\right):=\int_0^1 { t \log{(x-1)} \cdot (1-t)\log{(y-1)} \, dt}.
\label{eq:B_l1}
\end{equation}
But since a logarithmic function satisfies the functional equation
\begin{equation*}
l(x y)= l(x) + l(y),\qquad x,y \in (0,+\infty),
%\label{eq:}
\end{equation*}
and addition is the operation between $l(x)$ and $l(y)$ on the right hand side, arguing similarly when introducing the second additive pendant of the Beta function $\mathcal{A}_{2}$,
we should rather consider $\mathcal{L}_{2}: (1,+\infty) \to \R$ defined by
\begin{equation}
\mathcal{L}_{2} \left( x,y\right):=\int_0^1 { \left[t \log{(x-1)} + (1-t)\log{(y-1)} \right] \, dt}.
\label{eq:B_l2}
\end{equation}
as the logarithmic counterpart of the Euler Beta function.

It is not completely clear why in all these Beta functions the arguments are reduced by one, i.e. why $x-1$ and $y-1$, respectively, appear in the definiens of the Beta functions rather than $x$ and $y$ itself. 
A historic explanation could be that $\Gamma(n)=(n-1)!$ for a natural number $n$. 
Here as well we do not see any obvious reason
why the function $\Pi:(-1,+\infty) \to (0,+\infty)$ with $\Pi(x):=\Gamma{(x+1)}$ for $x>-1$, hence $\Pi(n)=n!$ for $n\in\N$, 
due to Gau{\ss} was not preferred as the natural extension of the sequence $(n!)_{n\in\N}$ by the majority of mathematicians%
\footnote{The domain of the function $\Pi$ is not symmetric with respect to zero, which one may consider as not so aesthetic. Moreover, nowadays the symbol $\Pi$ is usually avoided to denote functions since it is mostly used as a symbol for products. }.

%As the multiplicative and additive Beta functions, 
The two logarithmic Beta functions can be expressed by elementary functions since, for all $x,y \in (1, +\infty)$,
\begin{eqnarray*}
\mathcal{L}_{1} \left( x,y\right) &=& \log{(x-1) \log{(y-1) \int_0^1 {t } (1-t)} \, dt}\\
&=& \log{(x-1)} \log{(y-1)} \int_0^1  {{(t-t^2)} \, dt}\\
&=& \log{(x-1)} \log{(y-1)} {\left[\frac{t^2}{2}-\frac{t^3}{3}\right]}_0^1\\
&=& \frac{1}{6}\log{(x-1)} \log{(y-1)}.
%\label{eq:}
\end{eqnarray*}
Unsurprisingly when noticing the case of the function $\mathcal{A}_{2}$, we get
\begin{eqnarray*}
%\mathcal{B}_{a2} \left( x,y\right)=
\mathcal{L}_{2}\left( x,y\right)&=& \left[{\frac{t^2}{2} \log{(x-1)} + \left({t-\frac{t^2}{2}}\right)\log{(y-1)}}\right]_0^1\\
&=&{\frac{1^2}{2} \log{(x-1)} + \left({1-\frac{1^2}{2}}\right)\log{(y-1)}}\\
&=&\frac{1}{2} \log{(x-1)} + \frac{1}{2} \log{(y-1)}\\
&=&A(\log{(x-1)},\log{(y-1)}).
%\label{eq:B_l}
\end{eqnarray*}

%\iffalse
\section{Beta functions and their respective class of Cauchy quotients}
It is natural to ask whether these newly introduced Beta functions belong to the respective class of Cauchy quotients.
\subsection{The Euler Beta function - a beta-type function}
It is known that in the case of the Euler Beta function the answer is positive due to the relation between the Euler Beta and the Gamma function
 
\begin{equation}
\mathcal{B}\left( x,y\right) =
\frac{\Gamma(x) \Gamma(y)}{\Gamma(x+y)},\text{ \ \ \ \ \ }x,y \in (0, +\infty),
\label{eq:Beta_Gamma}
\end{equation}%
which means that Euler Beta function $\mathcal{B}$ is of the form 
\begin{equation*}
\frac{f(x) f(y)}{f(x+y)}
%\label{eq:}
\end{equation*}
for $f=\Gamma$ and thus belongs to the class of beta-type functions (cf. \cite{MatHim2}).

\subsection{The multiplicative Beta function - a multiplicative Cauchy quotient?}
The question whether the multiplicative Beta function $\mathcal{M}$ belongs to the class of multiplicative Cauchy quotients
can be reformulated as: is there a function $f:(1,+\infty) \to (0,+\infty)$ such that, for all $x,y \in (1, +\infty)$,

\begin{equation*}
\mathcal{M} \left( x,y\right)=\frac{f(x) f(y)}{f(xy)}, \qquad x \neq y,
%\label{eq:}
\end{equation*}
and \begin{equation*}
x-1=\frac{(f(x))^2}{f(x^2)} 
\end{equation*}
holds true?
\iffalse
From the last equation, since $f$ is a positive function, it follows that
\begin{equation*}
f(x)=\sqrt{(x-1)f(x^2)}, \qquad x \in (1, +\infty).
\end{equation*}
\fi

By \eqref{eq:B_m}, the multiplicative Beta function $\mathcal{M}$ coincides with the logarithmic mean (cf. \cite{Stolarsky}, \cite{Stolarsky2}), which is a Cauchy mean (cf. \cite{Matkowski}).
Hence, for all $x,y \in (1, +\infty)$,
\begin{equation*}
\frac{x-y}{\log{\left(\frac{x-1}{y-1}\right)}}=\frac{f(x) f(y)}{f(xy)}, \qquad x \neq y,
%\label{eq:}
\end{equation*}
and \begin{equation*}
x-1=\frac{(f(x))^2}{f(x^2)}.
\end{equation*}

\subsection{The additive Beta functions -  additive Cauchy quotients?}
The additive Beta function $\mathcal{A}_1$ (respectively $\mathcal{A}_2$) is an additive Cauchy quotient 
if there is an interval $I \subset \R$ which is closed under addition and a function $f: I \to (0, +\infty)$ such that, for all $x,y \in I$,
\begin{equation*}
\frac{1}{6}(x-1)(y-1)=\frac{f(x) + f(y)}{f(x+y)}, 
%\label{eq:}
\end{equation*}
(respectively, in case of  $\mathcal{A}_2$,
\begin{equation*}
\frac{(x-1)+(y-1)}{2}=\frac{f(x) + f(y)}{f(x+y)}, 
%\label{eq:}
\end{equation*}
)
holds true.
\iffalse
Setting $x=y$, it follows 
\begin{equation*}
f(x)=\frac{1}{12}(x-1) f(2x), \qquad x \in (0, +\infty)
%\label{eq:}
\end{equation*}
(respectively
\begin{equation*}
f(x)=\frac{x-1}{2} f(2x), \qquad x \in (0, +\infty)
%\label{eq:}
\end{equation*}).
\fi

\subsection{The logarithmic Beta function - a logarithmic Cauchy quotient?}
The logarithmic Beta function $\mathcal{L}_1$ (respectively $\mathcal{L}_2$) is a logarithmic Cauchy quotient
if there is an interval $I \subset (1,+\infty)$ which is closed under multiplication and a function $f: I \to (0, +\infty)$ such that, for all $x,y \in I$,
\begin{equation*}
\frac{1}{6}\log{(x-1)}\log{(y-1)}=\frac{f(x) + f(y)}{f(xy)}, 
%\label{eq:}
\end{equation*}
(respectively
\begin{equation*}
\frac{\log{(x-1)}+\log{(y-1)}}{2}=\frac{f(x) + f(y)}{f(xy)}, 
%\label{eq:}
\end{equation*}
)
holds true.
\iffalse
Setting $x=y$, it follows 
\begin{equation*}
f(x)=\frac{1}{12}\log{(x-1)} f(x^2), \qquad x \in I
%\label{eq:}
\end{equation*}
(respectively
\begin{equation*}
f(x)=\frac{\log(x-1)}{2} f(x^2), \qquad x \in I
%\label{eq:}
\end{equation*})
\fi

\section{Cauchy Beta functions: The case $k >2$}
In this section we introduce the Cauchy Beta functions of three or more variables.
Just as in the case of $\mathcal{M}$, the multiplicative Beta function of two variables,  
power functions having sum of exponents equal to one were multiplied with each other (since on the right hand side of the respective Cauchy equation $f(xy)=f(x) \cdot f(y)$ the numbers $f(x)$ and $f(y)$ are multiplied) and then integrated over the unit interval,
the multiplicative Beta function of three variables $\mathcal{M}_3:\left(1,+\infty \right)^{3}\rightarrow \R$ is
defined by

\begin{equation}
%\mathcal{B}_m\left( x,y\right) :=
\mathcal{M}_3\left( x,y,z\right)=
\int_0^1 \int_0^1 {(x-1)^s (y-1)^{t} (z-1)^{1-(s+t)}\, ds dt}, \text{ \ \ \ \ \ }x,y,z \in (1, +\infty).
\label{eq:Beta_m3:def}
\end{equation}

It clearly simplifies since, for all $x,y,z \in (1,+\infty)$, to
\begin{equation}
%\mathcal{B}_m\left( x,y\right) :=
\mathcal{M}_3\left( x,y,z\right)=
\frac{(x-z)(y-z)}{(z-1) \log{\left({\frac{x-1}{z-1}}\right)} \log{\left({\frac{y-1}{z-1}}\right)}}.
\end{equation}

Inductively, we obtain also the $k$-variable multiplicative Beta function
$\mathcal{M}_{k}:\left(1,+\infty \right)^{k}\rightarrow \R$, which for all $x_1,\ldots,x_k \in (1,+\infty)$ is
defined by
\begin{equation}
%\mathcal{B}_m\left( x,y\right) :=
\mathcal{M}_k\left( x_1, \ldots ,x_k\right)=
\int_0^1 \int_0^1 \cdots  \int_0^1 {(x_1-1)^{t_1} (x_2-1)^{t_2} \cdots (x_k-1)^{1-(t_1+ \cdots + t_{k-1})}
\, d{x_1} d{x_2} \ldots d{x_{k-1}}}.
\label{eq:Beta_mk:def}
\end{equation}

It can be written as
\begin{equation}
%\mathcal{B}_m\left( x,y\right) :=
\mathcal{M}_k\left( x_1, \ldots ,x_k\right)=
\frac{(x_1-x_k)(x_2-x_k)  \cdots (x_{k-1}-x_k) }{(x_k-1)^{k-1} \log{{\left(\frac{x_1-1}{x_k-1}\right)}} \cdot \cdots \cdot 
\log{\left({\frac{x_{k-1}-1}{x_k-1}}\right)}}
%\label{eq:Beta_mk:def}
\end{equation}

Analogously, the first kind three variable additive Beta function $\mathcal{A}_{1, 3}: \R^3 \to \R$ is defined by
\begin{equation}
%\mathcal{B}_m\left( x,y\right) :=
\mathcal{A}_{1, 3} \left( x, y ,z\right)=
\int_0^1 \int_0^1 {t(x-1)s(y-1)(1-(s+t))(z-1)}{\, dt ds},
%\label{eq:Beta_mk:def}
\end{equation}

simplifying to 
\begin{equation}
%\mathcal{B}_m\left( x,y\right) :=
\mathcal{A}_{1, 3} \left( x, y ,z\right)=
-\frac{1}{12} (x-1)(y-1)(z-1).
%\label{eq:Beta_mk:def}
\end{equation}

The first kind four variable additive Beta function $\mathcal{A}_{1, 4}: \R^4 \to \R$
simplifies to
\begin{equation}
%\mathcal{B}_m\left( x,y\right) :=
\mathcal{A}_{1, 4} \left( x, y ,z,w\right)=
-\frac{1}{8} (x-1)(y-1)(z-1)(w-1).
%\label{eq:Beta_mk:def}
\end{equation}

The first kind five variable additive Beta function $\mathcal{A}_{1, 5}: \R^5 \to \R$
simplifies to
\begin{equation}
%\mathcal{B}_m\left( x,y\right) :=
\mathcal{A}_{1, 5} \left( x, y ,z, w, a \right)=
-\frac{5}{48} (x-1)(y-1)(z-1)(w-1)(a-1).
%\label{eq:Beta_mk:def}
\end{equation}

\iffalse
%\begin{question}	
Can you find a closed formula for
the first kind $k$-variable additive Beta function $\mathcal{A}_{1, k}: \R^k \to \R$,
i.e. what is the number $c$ in front of the term $(x_1-1)(x_2-1) \cdots (x_k-1)$?
%\end{question}
\fi

\iffalse
The second (i.e. the more natural) additive Beta function of three variables $\mathcal{A}_{2, 3}: \R^3 \to \R$ defined by
\begin{equation}
\mathcal{A}_{2, 3} \left( x, y ,z\right)=
\int_0^1 {\int_0^1 {\left[t(x-1)+s(y-1)+(1-(s+t))(z-1)\right]}}{\, dt ds},
\end{equation}

simplifying to
\begin{equation}
\mathcal{A}_{2, 3} \left( x, y ,z\right)=
...
\end{equation}
What is general form of the second kind $k$-variable additive Beta form?
\fi

\iffalse
The two kinds of logarithmic Beta function of three variables $\mathcal{L}_{1, 3}, \mathcal{L}_{2, 3}: (1,+\infty)^k \to \R$ defined by
\begin{equation}
\mathcal{L}_{1, 3} \left( x, y ,z\right)=
\int_0^1 \int_0^1 {{\left[s\log(x-1) t\log(y-1)(1-(s+t))\log(z-1)\right]}}{\, dt ds},
\end{equation}
and
\begin{equation}
\mathcal{L}_{2, 3} \left( x, y ,z\right)=
\int_0^1 \int_0^1 {{\left[s\log(x-1) + t\log(y-1)+(1-(s+t))\log(z-1)\right]}}{\, dt ds},
\end{equation}
\fi

\section{Closing Remarks}

The questions whether the newly introduced pendants of Beta functions also
belong to their respective class of beta-type functions, 
namely whether $\mathcal{A}_{1,k} \in A_{f,k}$ (or, more importantly,  $\mathcal{A}_{2,k} \in A_{f,k}$)), 
$\mathcal{M}_{k} \in P_{f,k}$ and $\mathcal{L}_{1,k} \in L_{f,k}$ (or $\mathcal{L}_{2,k} \in L_{f,k}$)
will be answered in a future paper.

The Euler Beta function is the continuous analog of (the reciprocal of) a binomial coefficient
and consequently it appeas quite often in combinatorics.
By \eqref{eq:Beta_Gamma}, most properties of the Gamma function,
as, for instance, the duplication formula, Stirling's formula or its relation to the sine function, can be formulated
in terms of the Beta function.
Moreover, there are numerous applications of the Euler Beta function in applied sciences as String theory, Particle physics or Astrophysics.
As a result this note opens a can of new research questions, namely which properties of the Euler Beta function are shared 
by its here newly introduced siblings.
Also some problems of minor importance are left today without an answer (see sections $5.2$ to $5.4$ and $6$). 

From the perspective of functional equations also some probing issues strongly intrude themself, since one can introduce similar Beta functions not only for the Cauchy equations, but also related to some other functional equations.  
For example, let us consider the functional equation related to the sine-addition law
\begin{equation*}
	f(x+y)=f(x) g(y)+f(y)g(x).
\end{equation*}	  
Obviously,  $f=t \sin$, for $t \in \R$, and $g=\cos$ are continuous solutions on $\R$. 
The reader with an at least minor understanding of this paper will agree that 
it is more than natural to call $S: \R^2 \to \R$ defined by
\begin{equation*}
	S(x,y)=\int_0^1 \left({t \sin{x} \cos{y}+(1-t) \sin{y} \cos{x}}\right), \qquad x,y \in \R
\end{equation*}	  
(simplifying, for all $x,y \in \R$, to 
\begin{equation*}
	S(x,y)=\frac{1}{2} \sin{(x+y)},
\end{equation*})
the Beta function of the sine-addition functional equation (related to this given smooth solution pair).

We may summarize that the Euler Beta function indeed gives rise to introducing some kind of auxiliary Beta function 
in case of a functional euqation with solution involving some degree of freedom.

\end{document}